\documentclass[11pt,a4paper,reqno]{amsart}
\usepackage[T1]{fontenc}
\usepackage[utf8]{inputenc}
\usepackage{lmodern}
\usepackage[margin=28mm,headheight=14pt]{geometry}
\usepackage{amssymb,mathtools,booktabs,array,enumitem}
\usepackage[expansion=false]{microtype}
\usepackage{float,caption,tabularx}
\usepackage[unicode,hidelinks,pdfdisplaydoctitle]{hyperref}
\hypersetup{pdftitle={Infinitely many spectral gaps for Anderson--Bernoulli operators with sparse periodic backgrounds},pdfauthor={Cong Chen, Yong Li}}
\theoremstyle{plain}
\newtheorem{theorem}{Theorem}[section]
\newtheorem{existencetheorem}{Theorem}

\newtheorem{lemma}[theorem]{Lemma}
\newtheorem{proposition}[theorem]{Proposition}
\newtheorem{corollary}[theorem]{Corollary}
\theoremstyle{definition}
\newtheorem{definition}{Definition}[section]
\theoremstyle{remark}
\newtheorem{remark}[theorem]{Remark}

\DeclareMathOperator{\tr}{tr}
\DeclareMathOperator{\intr}{int}

\numberwithin{equation}{section}
\setlist{nosep,topsep=0.4em}
\allowdisplaybreaks[2]
\title[Spectral gaps for sparse periodic backgrounds]{Infinitely many spectral gaps for Anderson--Bernoulli operators with sparse periodic backgrounds}
\author{Cong Chen}
\thanks{Cong Chen: \href{mailto:chencong25@mails.jlu.edu.cn}{\texttt{chencong25@mails.jlu.edu.cn}}.}
\author{Yong Li}
\date{}
\begin{document}
\begin{abstract}
We consider discrete Anderson--Bernoulli operators with a periodic background. At a single coupling, every minimal period of at least three admits a background with at most four nonzero sites per period whose almost sure spectrum has infinitely many gaps. We give two constructions. The first uses at most seven nonzero sites per period and fixes the accumulation point at zero. For periods of at least five, the second uses exactly four nonzero sites per period and takes the lowest Dirichlet eigenvalue of the inserted chain as the accumulation point. In both constructions the gaps approach from the right. The widths of the selected gaps and their distances to the accumulation point are comparable to a common geometric factor divided by the period in the first family and by its cube in the second. Spectral intervals of comparable length separate successive gaps in each selected sequence.
\end{abstract}
\maketitle

\section{Introduction and main result}

A basic question in the spectral theory of Schr\"odinger operators is how the spatial structure of the potential determines the arrangement of bands and gaps. In one dimension, periodic potentials give finitely many spectral bands. Potentials generated by some aperiodic subshifts give Cantor spectra~\cite{DFII}. For an ergodic family the spectrum is almost surely a fixed set~\cite{DFI}, so one can ask how its gap structure changes when independent site disorder is added to a deterministic background.

Bellissard asked whether a bounded ergodic background with an Anderson perturbation necessarily has an almost sure spectrum with finitely many gaps. Avila, Damanik and Gorodetski proved finiteness for nondegenerate compactly supported noise when the background is sampled continuously along a homeomorphism of a connected compact base with an ergodic measure of full support~\cite{ADG}. A periodic background of minimal period \(p>1\) has a translation hull consisting of \(p\) discrete phases. The cyclic rotation of these phases together with the Bernoulli shift generates the random model considered here. For periodic backgrounds, Damanik and Gorodetski asked whether the spectrum of the perturbed operator is a finite union of intervals~\cite{DGInterval}. This question is stronger than asking whether the spectrum contains an interval. In the periodic background setting, the spectra of periodic noise realizations already supply intervals.

The first two background periods give finite unions of intervals. For a constant background, the Anderson spectrum formula gives at most two intervals in the Bernoulli case~\cite{DFII}. At period two, repeating the four binary words of length two gives the entire almost sure spectrum~\cite{DFG,WJMP}. At period three, a periodic noise word of length six can give an interval missed by all repetitions of words of length three~\cite{DFG}. Wood's dissertation gives an explicit example at period four and states the existence of infinitely many gaps for periods at least four~\cite{Wood}. Fillman and Wood construct a model of minimal period four whose gaps accumulate at zero from the left~\cite{FW}. Our earlier work establishes that three is the smallest possible background period for infinitely many gaps~\cite{CL}.

We ask whether such spectra can be realized at every prescribed minimal period \(p\ge3\), using one coupling and a bound on the number of nonzero background sites per period that is independent of \(p\). We also ask whether the selected gaps can accumulate at the same energy for all these periods. Sparse deterministic potentials can open all \(p-1\) bounded Floquet gaps with one nonzero site per period~\cite{Kaminaga}. In our setting, repeating a shorter background block leaves the minimal period unchanged, and every inserted site with zero background still carries an independent noise bit. The longer blocks must retain the connection that produces the gaps while accommodating every new noise word.

For a real periodic background \(v=(v_j)_{j\in\mathbb Z}\) and a coupling \(\lambda>0\), the total potential of the periodic Anderson--Bernoulli family~\cite{FW,CL} and the minimal period of its background are given by
\[
\begin{aligned}
 V_\omega(j)&=v_j+\lambda\omega_j\quad(j\in\mathbb Z),\qquad \omega\in\{0,1\}^{\mathbb Z},\\
 p&=\min\{q\in\mathbb Z_{>0}: v_{j+q}=v_j\text{ for all }j\in\mathbb Z\}.
\end{aligned}
\]
Setting \(\lambda=0\) gives the deterministic periodic model, and setting \(v\equiv0\) gives the Anderson--Bernoulli model. In our constructions, we prescribe \(p\) and keep the coupling \(\lambda\) fixed.

\noindent\textbf{Main result.}\enspace The discrete Schr\"odinger operator \(H_\omega\) with potential \(V_\omega\) acts on \(\ell^2(\mathbb Z)\) by
\begin{equation}
 (H_\omega u)(j)=u(j+1)+u(j-1)+(v_j+\lambda\omega_j)u(j).
 \label{eq:model}
\end{equation}
The variables \(\omega_j\) are independent and identically distributed, with \(\mathbb P(\omega_j=1)=\theta\) and \(\mathbb P(\omega_j=0)=1-\theta\), where \(0<\theta<1\). For fixed \(\theta\), the almost sure spectrum \(\Sigma(v)\) is deterministic and common to all phases of \(v\).

\par\addvspace{0.85\baselineskip}
\noindent\begin{minipage}{\linewidth}
\begin{existencetheorem}\label{thm:existence}
There is a coupling \(\lambda>0\) such that for every integer \(p\ge3\) one can choose a background \(v\) of minimal period \(p\) with
\[
 \#\{j\in\mathbb Z:0\le j<p,\ v_j\ne0\}\le4.
\]
For this background, \(\Sigma(v)\) has infinitely many gaps for every \(\theta\in(0,1)\).
\end{existencetheorem}

\begin{proof}
For \(p=3,4\), use the backgrounds in Theorem~\ref{thm:main}. For \(p\ge5\), Theorem~\ref{res:four} gives a background with exactly four nonzero sites per period at the same coupling.
\end{proof}
\end{minipage}
\par\addvspace{0.85\baselineskip}

Both constructions use the same coupling. The first gives gaps accumulating at zero from the right for every \(p\ge3\). The second uses exactly four nonzero background sites per period for \(p\ge5\), with a period-dependent accumulation energy. In the first construction, we use separate blocks for \(p=3,4,5\) and lengthen the zero run in one block for each parity when \(p\ge6\). Set
\begin{equation}
  \lambda=\frac{6083}{500},\qquad v_1=-\frac{11493}{1000},\qquad
  v_2=-\frac{153}{125},\qquad L=\frac{31}{2}.
  \label{eq:parameters}
\end{equation}
Let \(e\in(0.568030,0.568031)\) be the connection energy for \((0,v_1,v_2)\). At \(e\), the decaying solutions on two periodic tails join across a finite bridge, as specified in Lemma~\ref{lem:seed}. Set
\[
 \bar v_1=v_1-e,\qquad \bar v_2=v_2-e,\qquad h=1-L-e .
\]
Let \(0^{[m]}\) denote a string of \(m\) consecutive zeros, with \(0^{[0]}\) the empty string. Define the following blocks, taking \(k\ge0\) in the last two lines, by
\begin{equation}
\begin{aligned}
 V^{(3)}&=(-e,\bar v_1,\bar v_2),\\
 V^{(4)}&=(-1-e,-1,\bar v_1-1,\bar v_2),\\
 V^{(5)}&=(-L-e,0,L,\bar v_1,\bar v_2),\\
 V^{(6+2k)}&=(h,1,0^{[2k]},1,L,\bar v_1,\bar v_2),\\
 V^{(7+2k)}&=(h,2,1,0^{[2k]},2,L,\bar v_1,\bar v_2).
\end{aligned}
\label{eq:sparse}
\end{equation}

\begin{theorem}\label{thm:main}
Fix the coupling \(\lambda\) in~\eqref{eq:parameters}. For each integer \(p\ge3\), let \(v\) be the periodic extension of \(V^{(p)}\) in~\eqref{eq:sparse}, with \(v_n=V^{(p)}_{n\bmod p}\) and entries indexed from zero. The following two statements hold.
\begin{enumerate}[label=\textup{(\alph*)},leftmargin=*]
\item The minimal period is \(p\). For \(p\ge6\), the background takes values in \(\{0,1,2,L,h,\bar v_1,\bar v_2\}\) and has exactly six nonzero sites per period when \(p\) is even and seven when \(p\) is odd.
\item For every \(0<\theta<1\), the almost sure spectrum of~\eqref{eq:model} is the same set \(\Sigma_p=\Sigma(v)\) and contains zero. There are pairwise distinct bounded gaps \(G_{p,j}\) of \(\Sigma_p\) with
\[
 G_{p,j}\subset(0,\infty),\qquad \lim_{j\to\infty}\sup G_{p,j}=0 .
\]
\end{enumerate}
\end{theorem}

The proof of Theorem~\ref{thm:main} is given in Section~\ref{sec:transfer}.

To construct a background with four nonzero sites per period, we let the accumulation energy depend on \(p\). For \(p\ge5\), the lowest Dirichlet eigenvalue of the free chain on \(p-4\) sites is
\[
 E_p=-2\cos\!\left(\frac{\pi}{p-3}\right).
\]
The accumulation energy \(E_p\) tends to \(-2\) as \(p\to\infty\).

\begin{theorem}\label{res:four}
Use the coupling of Theorem~\ref{thm:main}. For each integer \(p\ge5\), extend the block
\begin{equation}\label{eq:four-block}
 \widetilde V^{(p)}=(-L-e,\,0^{[p-4]},\,L+E_p,\,\bar v_1+E_p,\,\bar v_2+E_p)
\end{equation}
periodically. The following two statements hold.
\begin{enumerate}[label=\textup{(\alph*)},leftmargin=*]
\item Its minimal period is \(p\), with exactly four nonzero background sites per period.
\item For every \(0<\theta<1\), the almost sure spectrum \(\widetilde\Sigma_p\) contains \(E_p\) and has infinitely many distinct bounded gaps accumulating at \(E_p\) from the right.
\end{enumerate}
\end{theorem}

The proof of Theorem~\ref{res:four} is given in Section~\ref{res:section}.

\begin{remark}\label{rem:theorem-a}
Theorem~\ref{thm:existence} is of interest for two reasons.
\begin{enumerate}[label=\textup{(\alph*)},leftmargin=*]
\item The fraction of nonzero background sites is at most \(4/p\). The entries in~\eqref{eq:four-block} are bounded uniformly in \(p\), so the mean absolute value of the background over a period is \(O(p^{-1})\). This gives examples with arbitrarily small background support density at the same coupling, for every nondegenerate Bernoulli law on \(\{0,1\}\).
\item To keep the connection energy at zero in Theorem~\ref{thm:main}, we insert pairs of sites with zero background into separate blocks for even and odd periods when \(p\ge6\). Allowing the energy to vary with \(p\) simplifies the construction. At the lowest Dirichlet eigenvalue \(E_p\), the transfer through the chain and its two flanking sites satisfies~\eqref{res:identity} when the noise vanishes at all those sites. A single block formula then works for both parities with four nonzero sites per period. If any noise bit on the chain is one, the endpoint Green function bound~\eqref{res:green} controls the outgoing slope uniformly in the chain length. For \(p\ge6\), the connection energy is \(E_p<0\). At \(p=5\), the chain has one site, so \(E_5=0\) and \(\widetilde V^{(5)}=V^{(5)}\). Both constructions then give the same background and connection.
\end{enumerate}
\end{remark}

\subsection*{Structure of the Proof}
\emph{Extending the period three connection.}
At the chosen energy \(E_*\), four block matrices agree with those of the period three model up to a common sign. For long periods, we obtain this agreement by inserting pairs of sites with zero background at \(E_*=0\), or by inserting a free chain and taking \(E_*=E_p\) to be its lowest Dirichlet eigenvalue. Strict inclusions in projective intervals control the remaining noise words.

\emph{Separating the gaps.}
For every sufficiently large \(n\), we match the two iteration scales. Multicone invariance gives a resolvent interval \(\mathcal J_{p,n}\), and a Floquet estimate for a periodic word gives an interval \(S_{p,n}\) in the spectral interior. The following order places the resolvent intervals in distinct bounded gaps.
\[
 \sup S_{p,n+1}<\inf\mathcal J_{p,n}
 \le\sup\mathcal J_{p,n}<\inf S_{p,n}.
\]

\emph{Measuring the gaps.}
The gap widths and their distances to \(E_*\) are of order \(R_B^{-2n}/d_p\). Here \(R_B^{-2n}\) comes from hyperbolic iteration and \(d_p\) is the energy derivative of the difference of the connecting slopes at \(E_*\). By the Wronskian identity, \(d_p\) equals the squared norm of the normalized connection solution. It is of order \(p\) in the zero energy family and \(p^3\) in the resonant family. The comparison constants are independent of \(p\), while the starting index may depend on the background. See Corollaries~\ref{cor:intro-quantitative} and~\ref{cor:resonant-gaps}.

\subsection*{The Organization of This Paper}
Section~\ref{sec:background} gives the spectral and projective preliminaries. Section~\ref{sec:seed} constructs the connection for period three and treats the short periods. Section~\ref{sec:sparse} extends the construction by inserting pairs of zero background sites. Section~\ref{sec:transfer} proves the gap separation and width estimates and applies them to the family with accumulation energy zero. Section~\ref{res:section} proves the resonant construction and its quantitative estimates. Appendix~\ref{app:width} verifies the numerical bounds for the integer scales. Appendix~\ref{app:algebraicity} proves that the connection energy of the period three model is algebraic.

\section{Spectral and projective preliminaries}\label{sec:background}

\subsection{Transfer}\label{pre:projective}

The total potential at site \(n\) is \(V_n=v_n+\lambda\omega_n\). For a solution \(u\) of the difference equation \(H_\omega u=Eu\), use the state \(X_n=(u(n),-u(n-1))^{\mathsf T}\). Transfer through a single site is given by
\[
 X_{n+1}=T_{E-V_n}X_n,\qquad
 T_a=\begin{pmatrix}a&1\\-1&0\end{pmatrix}.
\]
For $m<n$, the transfer from $m$ to $n$ is $T_{E-V_{n-1}}\cdots T_{E-V_m}$. Thus the first site in a spatial word contributes the rightmost factor.

We represent the line of a nonzero state \(X_n\) in \(\mathbb{RP}^1\) by its slope \(z=-u(n-1)/u(n)\), with \(z=\infty\) when \(u(n)=0\). For a matrix \(M\) of determinant one, \(Mz\) denotes the slope of the image line. For finite \(z\), let \(h_M(z)\) be the first coordinate of \(M(1,z)^{\mathsf T}\). Then
\[
 h_M(z)=M_{11}+M_{12}z,\qquad
 Mz=\frac{M_{21}+M_{22}z}{h_M(z)},\qquad
 \frac{d(Mz)}{dz}=h_M(z)^{-2}.
\]
In particular $T_a z=-1/(a+z)$. A zero of $h_M$ maps to infinity. The matrices $M$ and $-M$ have the same projective action.

A matrix $M\in\mathrm{SL}(2,\mathbb R)$ is hyperbolic when $|\operatorname{tr}M|>2$. Its unstable and stable lines, denoted by $u_M,s_M$, correspond to eigenvalues of moduli $\rho_M>1$ and $\rho_M^{-1}$, respectively. We use the same notation for their slopes in a finite chart. At a finite fixed slope $z$, the eigenvalue on $(1,z)^{\mathsf T}$ is $h_M(z)$. Expressing a vector in the eigenbasis gives
\begin{equation}\label{pre:ratio}
 \frac{Mz-u_M}{Mz-s_M}
       =\rho_M^{-2}\frac{z-u_M}{z-s_M}
\end{equation}
where both sides are defined. Each iteration of the block multiplies this ratio by the same scalar.

\begin{definition}\label{def:connection}
Let \(B(E)\) and \(A(E)\) be the transfer matrices over one period of the left and right periodic tails, respectively. Suppose both are hyperbolic near \(E_*\). Let \(D(E)\) be the transfer matrix through the finite block between the tails. A \emph{heteroclinic connection} at \(E_*\) is the equality \(D(E_*)u_{B(E_*)}=s_{A(E_*)}\) on the projective line. In a common local coordinate, let \(z_-(E)\) and \(z_+(E)\) denote \(D(E)u_{B(E)}\) and \(s_{A(E)}\). The connection is \emph{transverse} if \((z_--z_+)'(E_*)\ne0\), where the prime denotes differentiation with respect to \(E\).
\end{definition}

By the chain rule, changing coordinates multiplies the derivative of the difference by the nonzero derivative of the coordinate change at the common slope. At the connection, the common line gives a solution that decays at both ends. For the connection formulation, see \cite{FW}. The Wronskian proof of transversality is given in \cite{CL}.

For two solutions $u,v$ at the same energy, the Wronskian $u(n)v(n+1)-u(n+1)v(n)$ is constant in $n$, since it is the determinant of their states at $n+1$. For this calculation, let \(v=v(E)\) be a real solution that depends differentiably on \(E\). A dot denotes differentiation with respect to \(E\) at fixed total potential \(V\). Put
\begin{equation}\label{pre:wronskian}
 W_n=v(n)\dot v(n+1)-v(n+1)\dot v(n),\qquad
 W_n-W_{n-1}=v(n)^2 .
\end{equation}
Fix an energy interval on which both periodic tails are hyperbolic and the solutions that decay at the left and right ends \(v_-(n,E)\) and \(v_+(n,E)\) can be normalized to equal one at \(n=0\). Set \(z_\pm(E)=-v_\pm(-1,E)\).

On each decaying tail the solution is geometric at block boundaries, and differentiation introduces at most a linear factor. Thus \(W_n\) tends to zero at the corresponding end. Since \(W_{-1}=-\dot v(-1)\), summing~\eqref{pre:wronskian} on either side of the junction gives
\begin{equation}\label{pre:wr-sign}
 (z_--z_+)'(E)=\sum_{n<0}|v_-(n,E)|^2+\sum_{n\ge0}|v_+(n,E)|^2>0 .
\end{equation}
At a connection the two normalized solutions agree, so this derivative equals the squared norm in \(\ell^2\) of their common solution.

\subsection{Periodic spectral approximation}

A block consists of the sites \(kp,\ldots,kp+p-1\) for some \(k\in\mathbb Z\). It starts at background phase zero. Its noise word is \(\nu\in\{0,1\}^p\), and each such word has positive probability.

\begin{lemma}\label{pre:spectral}
The following two statements hold.
\begin{enumerate}[label=\textup{(\alph*)},leftmargin=*]
\item The spectrum $\Sigma$ of almost every realization is independent of \(\theta\in(0,1)\) and satisfies
\[
 \Sigma=\bigcup_{\omega}\sigma(H_\omega)
       =\overline{\bigcup_{\omega\ {\rm periodic}}\sigma(H_\omega)} .
\]
\item If the noise configuration has period \(q\), then \(\operatorname{lcm}(p,q)\) is a period of the total potential. It is a multiple of \(p\), so the \emph{monodromy matrix} \(P(E)\) through this common period is a product of whole block matrices. The trace \(\Delta(E)=\operatorname{tr}P(E)\) is the \emph{Floquet discriminant}. For real \(E\),
\[
 E\in\sigma(H_\omega)\quad\Longleftrightarrow\quad
             |\operatorname{tr}P(E)|\le2 .
\]
\end{enumerate}
\end{lemma}
\begin{proof}
For fixed \(\theta\in(0,1)\), the shift by \(p\) sites preserves \(v\) and is ergodic for the Bernoulli noise, so the almost sure spectrum is deterministic. Translation by one site conjugates the operators with adjacent background phases and preserves the noise law.

For the union over all configurations, use the Weyl sequence argument in~\cite{CL}, aligning the background phases by taking shifts in \(p\mathbb Z\) in place of \(3\mathbb Z\). Periodic approximation is~\cite{FW}, and the trace criterion is~\cite{DFII}.
\end{proof}

Each periodic spectrum is the closure of its interior, so periodic approximation implies that \(\operatorname{int}\Sigma\) is dense in \(\Sigma\)~\cite{FW}.

\subsection{Multicones and uniform hyperbolicity}\label{pre:multicones}

\begin{definition}\label{def:multicone}
A \emph{closed projective multicone} is a nonempty proper subset \(K\) of \(\mathbb{RP}^1\) consisting of finitely many pairwise disjoint nondegenerate closed arcs. For a finite family \(\mathcal F\subset\mathrm{SL}(2,\mathbb R)\), the multicone is \emph{forward invariant} if \(MK\subseteq K\) for every \(M\in\mathcal F\). It is \emph{strictly invariant under words of length \(q\)} if
\[
 M_q\cdots M_1K\subseteq\operatorname{int}K
 \qquad(M_1,\ldots,M_q\in\mathcal F).
\]
Here the interior is taken in the projective line. For \(q=1\) this is strict invariance under the generators.
\end{definition}

\begin{lemma}\label{pre:cone}
At a real energy $E_0$, suppose some multicone $K$ is strictly invariant for every product of $q$ block matrices, for a fixed $q\ge1$. Then a neighborhood of $E_0$ is disjoint from $\Sigma$.
\end{lemma}
\begin{proof}
Compactness and finiteness preserve all the strict inclusions on a common energy neighborhood. By the multicone criterion the family of products of length \(q\) is uniformly hyperbolic, and grouping the matrices in blocks of length \(q\) gives uniform hyperbolicity for the original family~\cite{ABY,FW}. The block form of Johnson's theorem excludes the neighborhood from \(\Sigma\)~\cite{FW}.
\end{proof}

\begin{remark}\label{rem:boundary-graph}
Let \(K\) be a forward invariant multicone with boundary \(\partial K\) in the projective line. The \emph{boundary transition graph} has vertex set \(\partial K\) and an edge \(x\to y\) whenever \(Mx=y\in\partial K\) for some \(M\in\mathcal F\). If this graph is acyclic, every word of length \(|\partial K|\), the number of boundary points, maps \(K\) into its interior~\cite{CL}.
\end{remark}

\section{The connection energy and projective intervals for period three}\label{sec:seed}

We locate the connection energy \(e\) for the background of period three and bound the images of the eight block maps. The insertions in Section~\ref{sec:sparse} preserve the four matrices needed at this energy up to a common sign. Section~\ref{sec:transfer} treats the energy dependence of the longer blocks.

We use the unshifted energy coordinate of the model of period three with background \((0,v_1,v_2)=(0,-11.493,-1.224)\). Terminating decimals denote exact rational numbers. For the three noise bits \(i,j,k\in\{0,1\}\) in spatial order, the transfer matrix in the notation of Section~\ref{pre:projective} is
\[
 M_{ijk}(E)=T_{E-v_2-\lambda k}T_{E-v_1-\lambda j}T_{E-\lambda i}.
\]
Set \(A=M_{000},R=M_{001},C=M_{010},B=M_{011},D=C^2\) and fix
\[
 w=-\tfrac32,\quad t=-\tfrac{307}{200},\quad x^-=-\tfrac14,\quad x^+=0,
 \qquad X=[x^-,0],\qquad\mathcal I=[w,0.104].
\]
The multicone begins at \(w\). The point \(t<w\) sets the left endpoint of the starting interval \([t,0]\) for the \(B\)-iterates. In expressions such as \(M0\), \(Bt\) and \(D0\), a matrix acts on the indicated slope.

Table~\ref{tab:seed-matrices} summarizes the roles of the four matrices.

\begin{table}[H]
\centering
\small
\begin{tabularx}{\linewidth}{@{}cc>{\raggedright\arraybackslash}X@{}}
\toprule
Matrix & Noise word & Role in the proof \\
\midrule
\(A\) & \(000\) & The stable slope \(s_A\) gives the tail decaying to the right. \\
\(R\) & \(001\) & The identity \(RA^{-1}C=B\) gives the boundary image \(R(A^{-1}C0)=B0\). \\
\(C\) & \(010\) & Two blocks with matrix \(C\) join the two slopes at \(E=e\), with \(C^2u_B=s_A\). \\
\(B\) & \(011\) & The unstable slope \(u_B\) gives the tail decaying to the left. \\
\bottomrule
\end{tabularx}
\caption{Transfer matrices for period three}
\label{tab:seed-matrices}
\end{table}

We use the energy dependent slopes \(\ell=Cx^-\), \(j_1^+=A^{-1}C0\), and \(\widehat i=u_B+\rho_B^{-4}\). The bounds in parts (b)--(d) below are at the connection energy \(e\).

\begin{lemma}\label{lem:seed}
The following statements hold.
\begin{enumerate}[label=\textup{(\alph*)},leftmargin=*]
\item There is a unique energy \(e\) in \(\mathcal J=[0.568030,0.568031]\) such that \(D u_B=s_A\).
\item At \(E=e\), the slopes satisfy
\begin{align}
 &s_B<-100<t<w<u_A<D0<s_A<\ell<x^-<0<u_B<\widehat i<0.104,\label{eq:seed-order}\\
 &0<Bt<0.08734,\qquad Bt<Bw<0.088<0.1<B0<u_B,\qquad D0>-0.53.\label{eq:seed-cut}
\end{align}
\item At \(E=e\), each projective map below has no pole on its stated source interval, and
\begin{align}
 C\mathcal I&\subset(w,x^-),&
 A[\ell,0.104]&\subset(w,D0),\nonumber\\
 R[w,j_1^+]&\subset(Bw,B0],&
 R[\ell,0.104]&\subset(Bw,B0).\label{eq:seed-returns}
\end{align}
In the first inclusion for \(R\), equality \(Rz=B0\) occurs only at \(z=j_1^+\).
\item At \(E=e\), on \([w,s_A]\), \([s_A,\infty)\), and \([t,\infty)\), respectively, the denominators satisfy \(h_A<0,h_{A^{-1}}<0,h_B>0\).
\end{enumerate}
\end{lemma}

\begin{proof}
The matrix entry and fixed slope bounds are given in Sections A.2 and A.3 of the appendix to~\cite{CL}. The entries are obtained from the product of three transfer matrices
\begin{equation}
 T_zT_yT_x=\begin{pmatrix}(zy-1)x-z&zy-1\\1-yx&-y\end{pmatrix}.
 \label{eq:seed-product}
\end{equation}
Throughout \(\mathcal J\) one has
\begin{gather}
 -0.55183<u_A<-0.55181,\quad -0.51438<s_A<-0.51436,\quad 0.10258<u_B<0.10260,\nonumber\\
 \rho_A^2<2.13,\qquad \rho_B>10.43,\qquad \widehat i<0.104 .\label{eq:seed-radii}
\end{gather}
The estimates for \(C\mathcal I\), \(A[\ell,0.104]\), and the two branches of \(R\), the bounds \(D0>-0.53\) and \(Bw<0.088<0.1<B0<u_B\), and the corresponding denominator signs are in Sections A.4 and A.6--A.8 of the appendix to~\cite{CL}.

To locate the connection, use the fixed point polynomial \(F_M(z)=M_{12}z^2+(M_{11}-M_{22})z-M_{21}\). For each energy in Table~\ref{tab:connection-slopes}, substitution into~\eqref{eq:seed-product} gives negative and positive values at the left and right endpoints, respectively, for \(F_A\) on the \(s_A\) interval and \(F_B\) on the \(u_B\) interval. The branches follow from \(-1<h_A<0\) and \(h_B>1\) throughout the respective intervals.

\begin{table}[H]
\centering
\small
\begin{tabular}{@{}>{$}c<{$}>{$}c<{$}>{$}c<{$}@{}}
\toprule
 E & s_A(E) & u_B(E)\\
\midrule
 0.568030 & (-0.5143698755,-0.5143698754) & (0.1025909043,0.1025909044)\\
 0.568031 & (-0.5143720045,-0.5143720044) & (0.1025909284,0.1025909285)\\
\bottomrule
\end{tabular}
\caption{Slope bounds at the endpoints of \(\mathcal J\)}
\label{tab:connection-slopes}
\end{table}

On these intervals the branch of \(D\) through \(u_B\) is increasing. Evaluating \(D\) at the two endpoints of the \(u_B\) interval and subtracting the opposite endpoints of the \(s_A\) interval gives the following bounds for \(\delta(E)=D(E)u_B(E)-s_A(E)\).
\[
\begin{aligned}
 \delta(0.568030)&\in(-1.38,-1.37)\times10^{-6},\\
 \delta(0.568031)&\in(1.38,1.39)\times10^{-6}.
\end{aligned}
\]
Throughout \(\mathcal J\), \(s_A\) and \(u_B\) are finite and continuous. So is \(Du_B\), because \(C\) has no pole on \(\mathcal I\) and maps it into itself. For the two \(C\) blocks between the \(B\) and \(A\) tails, \eqref{pre:wr-sign} gives \(\delta'(E)>0\) on \(\mathcal J\). Hence \(\delta\) has a unique zero \(e\) in the interior of \(\mathcal J\). At \(e\), monotonicity of \(D\) on \(\mathcal I\) gives \(D0<Du_B=s_A\).

At the new endpoint \(t=-307/200\), we have \(s_B<t\), \(B_{12}>0\), and \(h_B(s_B)>0\). Thus \(h_B>0\) throughout \([t,\infty)\) and \(Bt<Bw\). The formula \(Bt=(1-qv)/((cq-1)v-c)\) with \(q=e-0.673\), \(c=e-10.942\), and \(v=e+t\) gives
\[
 0<Bt<
 \frac{1-(0.104969)(0.966969)}{10.373969-(0.08896)(0.966970)}
 <0.08734 .
\]
\end{proof}

\subsection{The final two sites and the short periods}
The first \(p-2\) sites of a block of length \(p\) form its \emph{prefix}. A noise word is \emph{active} if all its prefix bits are zero. The last two bits are arbitrary, so there are four active words. The other words are \emph{inactive}. Section~\ref{sec:transfer} uses the three target intervals
\[
       (w,D0),\qquad (Bt,B0),\qquad (x^-,0).
\]

\begin{lemma}\label{lem:tail}
Suppose the slope entering the final two sites of the seed for period three is \(s\in(0.036,0.09)\). For noise symbols \((j,k)=(0,0),(0,1),(1,0),(1,1)\) at these sites, the outgoing slope lies in \((w,D0),(Bt,B0),(x^-,0),(x^-,0)\), respectively. Each inactive word of length three maps all of \(\mathcal I\) into the corresponding interval.
\end{lemma}
\begin{proof}
If \(j=0\), then \(12<e-v_1+s<14\), so the slope after the first step is \(v\in(-1/12,-1/14)\). If \(k=0\), then \(5/3<e-v_2+v<7/4\), giving an image in \((-3/5,-4/7)\). If \(k=1\), then \(10<-(e-v_2-\lambda+v)<100/9\), giving an image in \((0.09,0.1)\).

If \(j=1\), then \(-0.069<e-v_1-\lambda+s<0\) and the slope after the first step exceeds \(1000/69\). The denominator at the final step therefore exceeds \(1000/69-10.375>4\). The stated inclusions follow from~\eqref{eq:seed-cut}.

For any inactive word of length three, the noise symbol at the first site is \(1\). Its outgoing slope for input \(z\in\mathcal I\) is \((\lambda-e-z)^{-1}\in(1/14,1/11.4)\subset(0.036,0.09)\).
\end{proof}

\begin{lemma}\label{lem:four}
For the unshifted background of period four \((-1,e-1,v_1-1,v_2)\), the four active matrices at the critical energy \(e\) are \(A,R,C,B\), respectively. Every inactive matrix sends all of \(\mathcal I\) into one of the three open intervals above.
\end{lemma}
\begin{proof}
By the identity \(T_{y+1}T_1T_{x+1}=T_yT_x\), when the noise symbol at site \(1\) (counting from \(0\)) is zero, the product equals the original matrix for three sites determined by the remaining symbols. Lemma~\ref{lem:tail} covers the case in which the first symbol is \(1\).

If the symbol at the inserted site is \(1\), let \(i\) be the first noise bit and \(z\in\mathcal I\) the incoming slope. Set \(r=e+1+z-\lambda i\). After the first two sites, the slope is
\[
 U=\frac{r}{1+(\lambda-1)r},\qquad
 r\in[0.068030,1.672031]\ \text{or}\ [-12.097970,-10.493969].
\]
The denominator is positive on the first interval and negative on the second. The function increases on each interval. Substitution of the endpoints gives \(0.03866<U<0.091\). For example, at the lower end of the positive branch, \(0.068030>0.03866(1+11.166\cdot 0.068030)\). At the upper end of the negative branch use \(U=1/(11.166+1/r)\).

Let \(j,k\) be the last two noise bits and set \(Y=e-v_1+1-\lambda j+U,\ d=e-v_2-\lambda k\). The image after four steps is \(Z=-Y/(dY-1)\). The bounds for \(Y\) are
\[
\begin{aligned}
 13.09969<Y<13.152031 &\qquad (j=0),\\
 0.93369<Y<0.986031 &\qquad (j=1).
\end{aligned}
\]
If \(k=0\), then \(dY-1>1.792030(0.93369)-1>0\). If \(k=1\), then \(d<0<Y\). Since \(\partial_Y Z=(dY-1)^{-2},\ \partial_d Z=Y^2(dY-1)^{-2}\), it suffices to evaluate \(Z\) at the corners of each rectangle where both coordinates are minimal or both are maximal. Table~\ref{tab:period-four-slopes} lists the resulting bounds.
\begin{table}[H]
\centering
\small
\begin{tabular}{@{}>{$}c<{$}>{$}c<{$}>{$}c<{$}@{}}
\toprule
 j\backslash k &k=0&k=1\\ \midrule
 j=0&(-0.59,-0.58)&(0.095,0.097)\\
 j=1&(-1.40,-1.27)&(0.08737,0.089) \\
\bottomrule
\end{tabular}
\caption{Slope bounds for period four}
\label{tab:period-four-slopes}
\end{table}
The last entry has the smallest margin, since \(1/(10.373970+1/0.93369)>0.08737>0.08734>Bt\). The first column thus lies in \((w,D0)\) and the second in \((Bt,B0)\).
\end{proof}

\begin{lemma}\label{lem:five}
For the unshifted background of period five \((-L,e,e+L,v_1,v_2)\), the four active matrices at the critical energy \(e\) are \(-A,-R,-C,-B\), respectively. Every inactive word maps \(\mathcal I\) into the corresponding interval of Lemma~\ref{lem:tail}.
\end{lemma}
\begin{proof}
Let \(i,j\) be the noise bits at the two ends of the replacement block of length three and let \(z\in\mathcal I\) be its incoming slope. If the middle symbol is zero, direct multiplication gives
\[
  T_{-L-\lambda j}T_0T_{e+L-\lambda i}=-T_{e-\lambda(i+j)} .
\]
If \(i+j\in\{1,2\}\), the outgoing slope satisfies \(1/(\lambda(i+j)-e-z)\in(0.039,0.088)\), using \(w\le z\le 0.104\).

If the middle symbol is \(1\), set \(v=e+L-\lambda i+z>2.4\). The outgoing slope at the middle site is \(u=(\lambda+1/v)^{-1}\in(0,\lambda^{-1})\). At the right endpoint the outgoing slope lies between \((L+\lambda)^{-1}\) and \((L-\lambda^{-1})^{-1}\), both of which belong to \((0.036,0.09)\). Lemma~\ref{lem:tail} now applies.
\end{proof}

\section{Inserting pairs of sites with zero background}\label{sec:sparse}

Zero insertion also appears in the construction of spectra of Hausdorff dimension zero for generic limit-periodic potentials at all nonzero couplings~\cite{FGH}. Here the inserted sites carry Bernoulli noise.
Consider the backgrounds in~\eqref{eq:sparse} with \(p\ge6\). At \(E=0\), the transfer matrices at the last two sites \(T_{E-v_n-\lambda\omega_n}\) have parameters \(e-v_1-\lambda b\) and \(e-v_2-\lambda c\), respectively, where \(b,c\in\{0,1\}\). The two endpoints of the prefix have noiseless parameters
\[
 x=e+L-1,\qquad y=-L .
\]
The noiseless parameters in the interior of the prefix are
\[
 (-1,0^{[2k]},-1),\qquad (-2,-1,0^{[2k]},-2)
\]
for even and odd periods, respectively. At each site, a noise bit of \(1\) reduces the parameter by \(\lambda\).

\begin{lemma}\label{lem:pair}
Let \(E=0\) and suppose the incoming slope belongs to \(\mathcal I=[-3/2,\,0.104]\). The following two statements hold.
\begin{enumerate}[label=\textup{(\alph*)},leftmargin=*]
\item If all interior bits of the prefix vanish and its endpoint bits are \(i,j\), then the prefix matrix is \(\varepsilon_p T_{e-\lambda(i+j)}\), where \(\varepsilon_{6+2k}=(-1)^k\) and \(\varepsilon_{7+2k}=(-1)^{k+1}\).
\item If at least one bit anywhere in the prefix is \(1\), the outgoing slope lies in \((0.036,0.09)\), uniformly in \(k\).
\end{enumerate}
\end{lemma}
\begin{proof}
Let \(z\in\mathcal I\) be the incoming slope. Part (a) follows from the identities
\[
 M:=T_{-1}^2=\begin{pmatrix}0&-1\\1&-1\end{pmatrix},\qquad
 T_{-2}T_{-1}T_{-2}=-M,\qquad T_0^2=-I_2
\]
and \(T_{y-\lambda j} M T_{x-\lambda i}=T_{x+y+1-\lambda(i+j)}\). If all interior bits vanish and \(i+j\ge1\), the outgoing slope is \((\lambda(i+j)-e-z)^{-1}\). The estimate used in the proof of Lemma~\ref{lem:five} places this value in \((0.039,0.088)\).

Suppose now that at least one interior bit is \(1\). Set \(q=1\) for an even period and \(q=2\) for an odd period. Let \(\eta\) be the noise bit at the first interior site. At the left endpoint the denominator satisfies \(\tau=x-\lambda i+z\ge e+L-1-\lambda-3/2>1.4\), so the slope after the first interior site is
\[
 \xi=\frac{1}{q+\lambda\eta+1/\tau}\in(0,1/q).
\]
For an even period, \(\xi\) is the slope entering the run of pairs with zero background. For an odd period there is one more site, with noise bit \(\eta'\) and parameter \(-1-\lambda\eta'\), before the run, so the entering slope is \(\xi'=(1+\lambda\eta'-\xi)^{-1}\in(0,2)\).

A pair of sites with zero background and noise bits \(a,d\) acts on slopes by
\[
 F_{ad}(z)=\frac{z-\lambda a}{1+\lambda d(z-\lambda a)} .
\]
The pair acts by this formula even when the slope between its two sites is infinite. In particular, \(F_{00}\) is the identity. Put \(S=(-\infty,-8]\cup[0,9/100]\). The translation \(z\mapsto z-\lambda\) maps \((0,2)\) into \((-\infty,-8)\) and maps \(S\) into itself. The map \(z\mapsto z/(1+\lambda z)\) takes the nonnegative real axis into \([0,1/\lambda)\) and \((-\infty,-8]\) into \((1/\lambda,8/(8\lambda-1)]\), where \(8/(8\lambda-1)<9/100\). Thus the first pair whose bits differ from \(00\) sends the slope into \(S\), and every subsequent pair keeps it there.

Let \(z'\) be the slope entering the last interior site and \(\eta_{\rm last}\) its noise bit. The outgoing slope \(u\) then satisfies
\[
 0<u=\frac1{q+\lambda\eta_{\rm last}-z'}\le\frac{100}{91}<1.1,\qquad z'\in S.
\]

If every pair with zero background has bits \(00\), including the case \(k=0\), a noise bit of \(1\) must occur at one of the two or three remaining interior sites. If the bit at the last interior site is \(1\), its incoming slope is less than \(q\), so \(0<u<1/\lambda\).

Suppose that the last interior bit is \(0\). For an even period we have \(\xi<1/(1+\lambda)\), giving \(u=(1-\xi)^{-1}<1+1/\lambda\).

For an odd period, if the middle bit is \(1\), then \(\xi'<(\lambda+1/2)^{-1}<1\) and \(u=(2-\xi')^{-1}<1\). If the middle bit is \(0\), the first bit must be \(1\). In this case
\[
 \xi<\frac{1}{2+\lambda},\qquad \xi'<\frac{2+\lambda}{1+\lambda},
\]
again giving \(u<1+1/\lambda\).

In all cases \(u\in(0,1.1)\). At the right endpoint we obtain
\[
 0.036<\frac1{L+\lambda}<\frac1{L+\lambda j-u}<\frac1{L-1.1}<0.07 .
\]
\end{proof}

\begin{proposition}\label{prop:routes}
Fix \(p\ge3\) and use the background \(V^{(p)}\) in~\eqref{eq:sparse}. The following statements hold.
\begin{enumerate}[label=\textup{(\alph*)},leftmargin=*]
\item At energy zero, the four words with all prefix bits zero and final bits \(00,01,10,11\) give, respectively, the seed matrices \(A,R,C,B\) at energy \(e\), up to a common sign.
\item Set \(D=C^2\) and \(H=u_B+\rho_B^{-4}\) at the seed energy \(e\). The matrix of every other word maps \([w,H]\) compactly into one of \((w,D0)\), \((Bt,B0)\), and \(\intr X\). 
\item The four active matrix functions satisfy \(RA^{-1}C=B\) at every energy.
\end{enumerate}
\end{proposition}
\begin{proof}
At period three the four active words are those of the seed, and Lemma~\ref{lem:tail} controls the inactive words. Lemmas~\ref{lem:four} and~\ref{lem:five} treat the insertion of one site at period four and the single interior site at period five, respectively. For \(p\ge6\), Lemma~\ref{lem:pair} gives the asserted matrix identity when all prefix bits vanish. For every other word, the prefix maps \([w,H]\subset \mathcal I\) into \((0.036,0.09)\), and Lemma~\ref{lem:tail} sends this interval into one of the three open target intervals under the last two steps.

At energy \(E\), let \(M_{bc}\) be the active matrix with final bits \(b,c\). Let \(P\) be the common prefix matrix and \(V_b,U_c\) the transfer matrices at the penultimate and last sites, respectively. Then \(M_{bc}=U_cV_bP\), with any common sign absorbed into \(P\). Consequently,
\[
 M_{01}M_{00}^{-1}M_{10}
 =U_1V_0P\,P^{-1}V_0^{-1}U_0^{-1}\,U_0V_1P=M_{11}.\qedhere
\]
\end{proof}

\section{Period extension and separation of spectral gaps}\label{sec:transfer}

Finitely many interval inclusions at a critical energy suffice to prove the existence of infinitely many gaps. We use the hyperbolic iterates and periodic words from our earlier work~\cite{CL}. Extending the period changes the energy derivatives of the active matrices and introduces new words, since both noise values are allowed at every added site.

For the next proposition, let \(v\) be a real background of period \(p\), with almost sure spectrum \(\Sigma(v)\). For a binary noise word \(\nu\in\{0,1\}^p\), let \(M_\nu(E)\) be its transfer matrix over \(p\) sites. Designate four distinct words as active, with matrix functions \(A,R,C,B\). Use the coordinates of Section~\ref{pre:projective} and the constants \(t,w,x=x^-\) of Section~\ref{sec:seed}. Set \(D=C^2\) and \(H(E)=u_B(E)+\rho_B(E)^{-4}\). Here \(u_B\) is the unstable slope of \(B\), and \(\rho_B\) is its expanding eigenvalue modulus.

The index \(n\) labels the scales in the construction. Let \(\mathcal J_{p,n}\) denote the candidate resolvent interval defined by~\eqref{eq:transfer-loc}, \(G_{p,n}\) the gap that will contain it, and \(S_{p,n}\) the closed spectral interval at the same scale.

\begin{proposition}\label{prop:transfer}
Suppose an energy \(E_*\) satisfies the following conditions.
\begin{enumerate}[label=\textup{(\roman*)},leftmargin=*]
\item \textup{Seed matching.} At \(E_*\), the matrices \(A,R,C,B\) agree with the corresponding seed matrices of Lemma~\ref{lem:seed} at \(e\), up to a common sign.
\item \textup{Matrix identity.} The identity \(RA^{-1}C=B\) holds in a neighborhood of \(E_*\).
\item \textup{Inactive words.} For each remaining word \(\nu\), the image \(M_\nu(E_*)[w,H(E_*)]\) is compactly contained in one of
\[
 (w,D0),\qquad (Bt,B0),\qquad (x,0),
\]
where \(D0,Bt,B0\) are evaluated at \(E_*\).
\end{enumerate}

Let \(\delta_p(E)=D(E)u_B(E)-s_A(E)\) be the difference between the incoming unstable slope after the two \(C\) blocks and the stable slope of \(A\). Set \(R_B=\rho_B(E_*)\) and \(d_p=\delta'_p(E_*)\), with prime denoting the energy derivative. The following statements hold.
\begin{enumerate}[label=\textup{(\alph*)},leftmargin=*]
\item \(d_p>0\) and \(E_*\in\Sigma(v)\) for every \(0<\theta<1\).
\item There is a starting index \(N\), which may depend on the background, such that for every integer \(n\ge N\), \(\mathcal J_{p,n}\) lies in a bounded gap \(G_{p,n}\), and there is a closed interval \(S_{p,n}\subset\operatorname{int}\Sigma(v)\).
\item For \(n>N\), these intervals are ordered as
\[
 \sup G_{p,n}<\inf S_{p,n}\le\sup S_{p,n}<\inf G_{p,n-1}.
\]
\item For \(n\ge N\) and \(K\in\{G_{p,n},S_{p,n}\}\),
\[
\begin{aligned}
 |K|&\ge c_0\,\frac{R_B^{-2n}}{d_p},\\
 c_0\,\frac{R_B^{-2n}}{d_p}
 &\le \inf K-E_*\le\sup K-E_*\le C_0\,\frac{R_B^{-2n}}{d_p},
\end{aligned}
\]
where \(c_0,C_0>0\) depend only on the seed.
\item There is a constant \(\kappa>0\), depending only on the seed, such that every sufficiently small right neighborhood \((E_*,E_*+\varepsilon)\) contains a gap \(G_{p,n}\) and a spectral interval \(S_{p,n}\), each of length at least \(\kappa\varepsilon\). The threshold for \(\varepsilon\) may depend on the background.
\end{enumerate}
\end{proposition}

We keep the assumptions of Proposition~\ref{prop:transfer} throughout the following lemmas. First we normalize the active matrices and identify the energy derivative at the connection.

Let $\varepsilon\in\{1,-1\}$ be the common sign at $E_*$. Throughout a neighborhood of $E_*$ multiply each of the four physical block matrices by $\varepsilon$, so that the resulting $A,R,C,B$ coincide with the seed at $E_*$. This normalization preserves the identity $RA^{-1}C=B$ and the projective maps. For a word of $q$ active blocks it multiplies the trace by $\varepsilon^q$, preserving its absolute value.

For \(M=A,B\), let \(u_M,s_M\) be the unstable and stable slopes, and let \(\rho_M\) denote the modulus of the unstable eigenvalue. These quantities are real analytic near \(E_*\).

\begin{lemma}\label{lem:transfer-transverse}
The connection is transverse. More precisely, \(\delta_p(E_*)=0\), and \(d_p=\delta'_p(E_*)\) equals the squared \(\ell^2\)-norm of the connection solution \(\psi\) normalized by \(\psi(0)=1\) at the start of the right \(A\) tail. In particular, \(d_p>0\).
\end{lemma}

\begin{proof}
The seed connection gives $\delta_p(E_*)=0$. Form a potential with a periodic $B$ tail on the left, two $C$ blocks, and a periodic $A$ tail on the right. Set the origin at the start of the first $A$ block. Normalize the solutions decaying at $-\infty$ and $+\infty$ to equal $1$ at the origin. Their slopes there are $Du_B$ and $s_A$, respectively. At $E_*$ the slopes agree, so the two solutions join to give a real solution $\psi\in\ell^2(\mathbb Z)$.

Applying the Wronskian identity~\eqref{pre:wr-sign} at the junction gives
\begin{equation}\label{eq:transfer-wr}
 d_p:=\delta'_p(E_*)=\sum_{j\in\mathbb Z}|\psi(j)|^2>0 .
\end{equation}
Thus $\delta_p$ is a local coordinate to the right of $E_*$.
\end{proof}

The derivative converts projective displacements into energy increments. We next match the two integer scales and locate an energy interval on which the endpoints have the required order.

A superscript $*$ denotes evaluation at $E_*$. Set $R_A=\rho_A^*$, $R_B=\rho_B^*$, $u=u_B^*$, $\sigma=s_B^*$, $\Delta=s_A^*-u_A^*$, and $d_0=\partial_zD^*(u)>0$. The coefficients for the two families of directions are
\[
\begin{aligned}
 P(z)&=\Delta\frac{z-s_A^*}{z-u_A^*},&
 P_-&=P(C^*x),& P_+&=P(C^*0),\\
 Q_-&=d_0(u-\sigma)\frac{u-t}{t-\sigma},&
 Q_+&=R_B^2d_0(u-\sigma)\frac{u}{-\sigma}.
\end{aligned}
\]
Take \(r_0,r_1,\eta,T\) as in Lemma~\ref{lem:width-parameters}.
In the periodic word \(B^nC^2A^{m_n}\), \(n\) counts the \(B\) blocks and \(m_n\) counts the \(A\) blocks. For sufficiently large \(n\), define
\begin{equation}\label{eq:transfer-integers}
 m=m_n=\min\{j\in\mathbb Z_{\ge0}:R_A^{2j}R_B^{-2n}\ge r_0\},\quad
 \alpha_n=R_A^{-2m},\quad\beta_n=R_B^{-2n},\quad r_n=\beta_n/\alpha_n .
\end{equation}

\begin{lemma}\label{lem:transfer-window}
For all sufficiently large \(n\), the integers and scales in~\eqref{eq:transfer-integers} satisfy \(r_0\le r_n<R_A^2r_0\le r_1\), \(m\asymp n\) and \(n\alpha_n\to0\). The interval \(\mathcal J_{p,n}\) of~\eqref{eq:transfer-loc} satisfies
\[
 E_*+\frac{P_-\alpha_n}{d_p}<\inf\mathcal J_{p,n}
 \le\sup\mathcal J_{p,n}<E_*+\frac{P_+\alpha_n}{d_p}.
\]
For every \(E\in\mathcal J_{p,n}\), the endpoint order~\eqref{eq:transfer-order} holds.
\end{lemma}

\begin{proof}
The minimality of $m$ gives $r_0\le r_n<R_A^2r_0\le r_1$, so $m\asymp n$ and $n\alpha_n\to0$. Let
\[
 g(r)=\tfrac12\bigl(P_-+Q_-r+\min\{Q_+r,P_+\}\bigr).
\]
By Lemma~\ref{lem:width-parameters}, \(g(r)\) is more than \(3\eta/2\) from both endpoints of \((P_-+Q_-r,\min\{P_+,Q_+r\})\) for all \(r\in[r_0,r_1]\).

For $|E-E_*|=O_p(\alpha_n)$, analyticity and the bound $m=O(n)$ give
\[
 \rho_A(E)^{-2m}/\alpha_n=1+O_p(n\alpha_n),\qquad
 \rho_B(E)^{-2n}/\beta_n=1+O_p(n\alpha_n).
\]
The subscript $p$ allows the implicit constant to depend on the fixed periodic background. At energy $E$, set $j_m^-=A^{-m}Cx$, $j_m^+=A^{-m}C0$, and $h_n=u_B+\rho_B^{-4n}$. Apply the identity for the cross ratio in~\eqref{pre:ratio} to \(B^n\) and to \(A^{-m}\), interchanging the two fixed points for \(A^{-m}\). Expanding \(D\) at \(u_B\) gives
\begin{equation}\label{eq:transfer-asym}
\begin{aligned}
 j_m^\pm-s_A&=P_\pm\alpha_n+O_p(n\alpha_n^2),&
 DB^{n-1}0-s_A&=\delta_p-Q_+\beta_n+O_p(n\alpha_n^2),\\
 DB^nt-s_A&=\delta_p-Q_-\beta_n+O_p(n\alpha_n^2),&
 Dh_n-s_A&=\delta_p+O_p(\alpha_n^2).
\end{aligned}
\end{equation}
For example, $B^nt-u_B=-(u_B-s_B)(u_B-t)(t-s_B)^{-1}\rho_B^{-2n}+O_p(\rho_B^{-4n})$. The exponent $n-1$ accounts for the factor $R_B^2$ in $Q_+$. All four expansions are uniform in the stated energy range.

Consider the closed interval
\begin{equation}\label{eq:transfer-loc}
 \mathcal J_{p,n}=E_*+\frac{\alpha_n}{d_p}
     [g(r_n)-\eta,\ g(r_n)+\eta].
\end{equation}
For $E\in\mathcal J_{p,n}$, we have $\delta_p(E)/\alpha_n=d_p(E-E_*)/\alpha_n+O_p(\alpha_n)$. Substitution into \eqref{eq:transfer-asym} gives, for all sufficiently large $n$,
\begin{equation}\label{eq:transfer-order}
 DB^{n-1}0<s_A<j_m^-<DB^nt<Dh_n<j_m^+ .
\end{equation}
The inequality $DB^nt<Dh_n$ follows directly from $B^nt<u_B<h_n$ and the strict monotonicity of $D$.

Since \(Q_->0\) and \(r_n\ge r_0>0\), the margin at the two endpoints of the matching interval gives \(P_-<g(r_n)-\eta\) and \(g(r_n)+\eta<P_+\). Substituting these bounds into~\eqref{eq:transfer-loc} proves the stated containment.
\end{proof}

The endpoint inequalities in Lemma~\ref{lem:transfer-window} determine the placement of the components of the multicone. Take \(n\) sufficiently large that \(m\ge2\) and \(n\ge4\). Fix \(E\in\mathcal J_{p,n}\) and evaluate the matrices and slopes below at \(E\). Set \(X=[x,0]\) and define
\[
\begin{gathered}
 W=[w,DB^{n-1}0],\qquad L_0=[Cx,CB0],\qquad J_j=A^{-j}CX\ (1\le j\le m),\\
 I_i=B^i[t,0]\ (1\le i<n),\qquad I_n^*=[B^nt,h_n],\\
 Y_i=CI_i\ (2\le i<n),\qquad Y_n^*=CI_n^* .
\end{gathered}
\]
In $I_n^*$ and $Y_n^*$ the star marks the final intervals, whose right endpoints have been extended. Let $K$ be the union of these intervals and $X$. Let \(j_j^\pm\) be the endpoints of \(J_j\).

\begin{lemma}\label{lem:transfer-cone-order}
For all sufficiently large \(n\) and every \(E\in\mathcal J_{p,n}\), the set \(K\) is a union of pairwise disjoint nondegenerate closed intervals in \([w,H]\). Their order on the real line is
\[
 W,\ J_m,\ldots,J_1,\ L_0,\ Y_2,\ldots,Y_{n-1},Y_n^*,\ X,\ I_1,\ldots,I_{n-1},I_n^* .
\]
\end{lemma}

\begin{proof}
The strict inclusions and inequalities of Lemma~\ref{lem:seed} persist for \(E\) near \(E_*\). $B$ is increasing on $[t,H]$ and $t<0<Bt<B0<u_B$. Iterating gives $B^i0<B^{i+1}t$. Applying $C$ separates the $Y_i$ and gives $CB0<CB^2t$. Since $Ch_n\le CH<x$, the interval $Y_n^*$ lies to the left of $X$.

To order the \(J_j\), set \(\ell=Cx\). By Lemma~\ref{lem:seed}, \(\ell>s_A\) and \(A\ell<D0<s_A=As_A\). Since \(A\) is increasing on each branch, it has a pole \(z_A\in(s_A,\ell)\). Its affine denominator therefore satisfies \(h_A(s_A)<0<h_A(\ell)\), with \(A_{12}>0\). The inverse branch maps \([s_A,\infty)\) increasingly onto \([s_A,z_A)\). Thus \(s_A<j_1^-<j_1^+<\ell\). Applying \(A^{-j}\) to \(j_1^+<\ell\) gives \(j_{j+1}^+<j_j^-\) for \(1\le j<m\). Finally, \eqref{eq:transfer-order} separates $W$ from $J_m$. Thus $K$ is a finite union of pairwise disjoint nondegenerate closed intervals in $[w,H]$.
\end{proof}

With the components ordered by Lemma~\ref{lem:transfer-cone-order}, we track their endpoints under the active matrices.

\begin{lemma}\label{lem:transfer-resolvent}
For all sufficiently large \(n\) and every \(E\in\mathcal J_{p,n}\), the multicone \(K\) is invariant under the block matrices. There is a length \(q\) for which every word of length \(q\) maps \(K\) into its interior. In particular, \(\mathcal J_{p,n}\cap\Sigma(v)=\varnothing\).
\end{lemma}

\begin{proof}
Let $i_j^\pm,j_j^\pm,y_j^\pm$ denote the endpoints of the corresponding intervals, $\ell^\pm$ those of $L_0$, and $w_n^+$ the right endpoint of $W$. Hats denote the endpoints of $I_n^*$ and $Y_n^*$. With the convention of Remark~\ref{rem:boundary-graph}, Table~\ref{tab:multicone-transitions} lists all transitions from $\partial K$ to $\partial K$ under the four active matrices.
\begin{table}[htbp]
\centering
\renewcommand{\arraystretch}{1.14}
\begin{tabularx}{\linewidth}{@{}>{\centering\arraybackslash}p{.12\linewidth}@{\hspace{1.4em}}>{\raggedright\arraybackslash}X@{}}
\toprule
Matrix & Boundary transitions \\
\midrule
\(A\) & \(j_j^\pm\mapsto j_{j-1}^\pm\ (2\le j\le m),\quad j_1^-\mapsto\ell^-\) \\
\(R\) & \(j_1^+\mapsto i_1^+\) \\
\(B\) & \(0\mapsto i_1^+,\quad i_j^\pm\mapsto i_{j+1}^\pm\ (1\le j\le n-2),\quad i_{n-1}^-\mapsto\widehat i_n^-\) \\
\(C\) & \(x\mapsto\ell^-,\quad i_1^+\mapsto\ell^+,\quad i_j^\pm\mapsto y_j^\pm\ (2\le j<n),\quad \widehat i_n^\pm\mapsto\widehat y_n^\pm,\qquad y_{n-1}^+\mapsto w_n^+\) \\
\bottomrule
\end{tabularx}
\caption{Boundary transitions of \(K\)}
\label{tab:multicone-transitions}
\end{table}

For \(M=A,B\), factoring the numerator at the two fixed slopes gives
\[
 Mz-z=-\frac{M_{12}(z-u_M)(z-s_M)}{h_M(z)} .
\]

On \(W\) we have \(w<u_A<D0<w_n^+<s_A\) and \(h_A<0\). The formula gives \(Aw>w\) and \(Aw_n^+<w_n^+\), so monotonicity yields \(AW\subset\operatorname{int}W\). On the \(J\) intervals,
\[
 AJ_j=J_{j-1}\quad (2\le j\le m),\qquad AJ_1=CX\subset L_0 .
\]
Here \(C0<CB0\), so only the left endpoint of \(CX\) meets the boundary of \(L_0\). Every remaining component lies in \([Cx,H]\) and has image in \((w,D0)\).

The components outside the \(I\) intervals lie in \([w,0]\). Their images under \(B\) are in \((Bt,B0]\), since \(Bt<Bw\), with \(0\) mapping to \(B0\). The ordinary intervals satisfy
\[
 BI_j=I_{j+1}\quad (1\le j\le n-2),\qquad BI_{n-1}\subset I_n^* .
\]
Only the left endpoint of \(BI_{n-1}\) meets the boundary of \(I_n^*\). On \(I_n^*\) we have \(s_B<B^nt<u_B<h_n\), \(B_{12}>0\), and \(h_B>0\). The fixed slope formula gives \(B(B^nt)>B^nt\) and \(Bh_n<h_n\), so \(BI_n^*\subset\operatorname{int} I_n^*\).

Since \(x<0<Bt<B0\), the map \(C\) sends both \(X\) and \(I_1\) into \(L_0\). From \(Cw>w\) and \(j_1^+<Cx<C0\) we obtain \(C[w,j_1^+]\subset(w,D0)\). Also \(CL_0=D[x,B0]\subset(w,w_n^+)\), by \(D[w,H]\subset(w,x)\) and \(B0<B^{n-1}0\).

For each ordinary \(Y_j\), the image \(CY_j=DI_j\) lies in \((w,w_n^+]\) and meets \(w_n^+\) only for \(j=n-1\). The final image satisfies \(CY_n^*=DI_n^*\subset\operatorname{int}J_m\) by~\eqref{eq:transfer-order}.

For \(R\), use the inclusion \(K\subset[w,j_1^+]\cup[Cx,H]\). On each of these intervals \(R\) has no pole. Its endpoint values satisfy
\[
 Rw,\ R(Cx)>Bw>Bt,\qquad RH<B0,\qquad Rj_1^+=B0 .
\]
The last equality follows from \(RA^{-1}C=B\). Thus \(RK\subset I_1\), with the sole boundary contact listed in the table.

For an inactive matrix \(M\), parameterize the source interval by \(z=w+s(H(E)-w)\), \(0\le s\le1\). At \(E_*\), \(|h_M(z)|\) and the distance of \(Mz\) from either endpoint of its target interval have positive lower bounds. The target endpoints vary continuously with \(E\). By compactness in \(s\), the same inclusion holds near \(E_*\). For fixed \(p\) there are only finitely many inactive words, so one neighborhood suffices for all of them.

The inclusions established for \([w,H(E)]\) apply to \(K\subset[w,H(E)]\). As \(D0<w_n^+\), the three target intervals \((w,D0)\), \((Bt,B0)\), and \((x,0)\) lie in the interiors of \(W\), \(I_1\), and \(X\), respectively. Hence each inactive word maps \(K\) into its interior.

In Table~\ref{tab:multicone-transitions}, a boundary path through the \(J\) intervals has decreasing index and can exit only at an endpoint of \(L_0\) or \(I_1\). Through the \(I\) intervals the index increases, with exits to a \(Y\) interval or \(L_0\). From a \(Y\) endpoint there is at most one further edge, ending at a vertex with no outgoing edge. Thus the graph is acyclic, as in the calculation for period three in~\cite{CL}.

Remark~\ref{rem:boundary-graph} now gives strict invariance under all words of a fixed length. Lemma~\ref{pre:cone} applies at every \(E\in\mathcal J_{p,n}\), so \(\mathcal J_{p,n}\cap\Sigma(v)=\varnothing\).
\end{proof}

Lemma~\ref{lem:transfer-resolvent} excludes \(\mathcal J_{p,n}\) from the spectrum. To distinguish the gaps at successive scales, we now locate spectral intervals using the Floquet discriminant.

\begin{lemma}\label{lem:transfer-spectral}
For every sufficiently large \(n\), there is a closed interval \(S_{p,n}\subset\operatorname{int}\Sigma(v)\) with
\[
 \sup\mathcal J_{p,n}<\inf S_{p,n}<\sup S_{p,n}
 \le E_*+\frac{2T\alpha_n}{d_p},\qquad
 |S_{p,n}|\ge\frac{\xi\alpha_n}{d_p},
\]
where \(\xi>0\) depends only on the seed. Successive scales satisfy
\(\sup S_{p,n+1}<\inf\mathcal J_{p,n}\).
\end{lemma}

\begin{proof}
Repeat the spatial word $B^nC^2A^m$ periodically, reading from left to right. Its transfer matrix over one period agrees with $A^{m_n}DB^n$ up to sign. As in~\cite{CL}, let \(\Pi_M^+,\Pi_M^-\) be the projections onto the unstable and stable directions of \(M=A,B\). Put \(\tau_{ij}=\operatorname{tr}(\Pi_A^i D\Pi_B^j)\) for \(i,j\in\{+,-\}\). Set \(v(z)=(1,z)^{\mathsf T}\). The projection formula \(\Pi_M^+v(z)=\frac{z-s_M}{u_M-s_M}v(u_M)\) and \(Dv(z)=h_D(z)v(Dz)\) give
\[
 \tau_{++}=-H_0\delta_p,\qquad
 H_0=\frac{h_D(u_B)(s_B-u_A)}
 {(s_A-u_A)(s_B-u_B)}\ne0 .
\]
The normalized matrices \(A\) and \(B\) have negative and positive eigenvalues, respectively, since \(h_A(s_A)<0<h_B(u_B)\). Expand \(A^m\) and \(B^n\) in their spectral projections. After division by the nonzero factor $(-1)^m\rho_A^m\rho_B^nH_0$, the trace is exactly
\[
 -\delta_p+\rho_A^{-2m}a+\rho_B^{-2n}b+\rho_A^{-2m}\rho_B^{-2n}c,
 \qquad (a,b,c)=(\tau_{-+},\tau_{+-},\tau_{--})/H_0 .
\]
At \(E_*\) we have
\[
 D^*v(u)=h_{D^*}(u)v(s_A^*),\qquad
 \det(D^*v(u),D^*v(\sigma))=\sigma-u .
\]
The first identity gives \(\tau_{-+}^*=h_{D^*}(u)(s_A^*-\sigma)/(u-\sigma)\). The determinant gives the coefficient of \(v(u_A^*)\) in \(D^*v(\sigma)\) in the eigenbasis of \(A^*\), namely \((u-\sigma)/(\Delta h_{D^*}(u))\). Hence \(\tau_{+-}^*=(u-u_A^*)/(\Delta h_{D^*}(u))\). Dividing both coefficients by \(H_0^*\) gives
\[
 a^*=\Delta+\frac{\Delta^2}{u_A^*-\sigma}>\Delta>P_+,
 \qquad b^*=\frac{(u-u_A^*)(u-\sigma)}
 {h_{D^*}(u)^2(u_A^*-\sigma)}>0 .
\]
Lemma~\ref{lem:width-parameters} gives \(T>a^*+r_1b^*\). When \(\delta_p(E)=s\alpha_n\) with \(0\le s\le T\), division by \(\alpha_n\) gives
\begin{equation}\label{eq:transfer-trace}
 \frac{\operatorname{tr}(A^{m_n}DB^n)}
 {(-1)^m\rho_A^m\rho_B^nH_0\alpha_n}
   =-s+a^*+r_nb^*+O_p(n\alpha_n)
\end{equation}
uniformly in \(s\). Set \(s_n=a^*+r_nb^*\). By~\eqref{eq:transfer-integers}, \(R_A^{m_n}R_B^n\alpha_n=r_n^{-1/2}\). The estimates preceding~\eqref{eq:transfer-asym} and continuity of \(H_0\) thus give
\[
 |\rho_A^{m_n}\rho_B^n H_0\alpha_n|\le M,\qquad
 M=\frac{2|H_0^*|}{\sqrt{r_0}},
\]
for all large \(n\), uniformly on \(0\le s\le T\). Choose a seed constant \(\xi>0\) smaller than \(a^*-P_+\), \(T-a^*-r_1b^*\), and \(1/(2M)\). Let \(S_{p,n}\) be the image of \(\alpha_n[s_n-\xi,s_n+\xi]\) under the local inverse of \(\delta_p\) at \(E_*\). On this closed interval, \eqref{eq:transfer-trace} bounds the absolute trace by \(M(\xi+O_p(n\alpha_n))<2\) for large \(n\). By continuity and Lemma~\ref{pre:spectral}, \(S_{p,n}\) lies in the interior of the spectrum of the periodic word and hence in the interior of \(\Sigma(v)\).

The inverse has derivative \(1/d_p\) at zero. Consequently \(|S_{p,n}|=(2\xi+o_p(1))\alpha_n/d_p\), where the limit is taken with the background fixed. Thus \(|S_{p,n}|\ge\xi\alpha_n/d_p\) for all sufficiently large \(n\). Since \(g(r_n)+\eta<P_+<s_n-\xi<s_n+\xi<T\) with a uniform margin at \(P_+\), the expansion \(E=E_*+\delta_p(E)/d_p+O_p(\delta_p(E)^2)\) gives
\begin{equation}\label{eq:transfer-separator}
\begin{gathered}
 \mathcal J_{p,n}\cap\Sigma(v)=\varnothing,\qquad S_{p,n}\subset\operatorname{int}\Sigma(v),\\
 \sup\mathcal J_{p,n}<\inf S_{p,n}<\sup S_{p,n}\le E_*+2T\alpha_n/d_p .
\end{gathered}
\end{equation}
By~\eqref{eq:transfer-integers} and Lemma~\ref{lem:width-parameters},
\[
 \frac{\alpha_{n+1}}{\alpha_n}=\frac{r_n}{r_{n+1}R_B^2}
 \le \frac{r_1}{r_0R_B^2}<0.02,\qquad
 2T\cdot 0.02=0.008<P_-<g(r_n)-\eta .
\]
By~\eqref{eq:transfer-loc} and~\eqref{eq:transfer-separator}, \(\sup S_{p,n+1}<\inf\mathcal J_{p,n}\).
\end{proof}

\begin{proof}[\normalfont\bfseries Proof of Proposition~\ref{prop:transfer}]
Lemma~\ref{lem:transfer-transverse} gives \(d_p>0\). The spectral intervals of Lemma~\ref{lem:transfer-spectral} approach \(E_*\), so \(E_*\in\Sigma(v)\) by closedness.

For all sufficiently large \(n\), Lemmas~\ref{lem:transfer-resolvent} and~\ref{lem:transfer-spectral} place the resolvent interval \(\mathcal J_{p,n}\) between the two spectral intervals \(S_{p,n+1}\) and \(S_{p,n}\). Since both spectral intervals lie in the interior of \(\Sigma(v)\), the component \(G_{p,n}\) containing \(\mathcal J_{p,n}\) is bounded and satisfies
\[
 \sup S_{p,n+1}<\inf G_{p,n}<\sup G_{p,n}<\inf S_{p,n}.
\]
In particular, \(S_{p,n}\) separates \(G_{p,n}\) from \(G_{p,n-1}\), and all these gaps are distinct.

By Lemma~\ref{lem:transfer-window} and~\eqref{eq:transfer-separator}, both intervals \(K\in\{G_{p,n},S_{p,n}\}\) satisfy
\[
 E_*+\frac{P_-\alpha_{n+1}}{d_p}<\inf K\le\sup K\le E_*+\frac{2T\alpha_n}{d_p}.
\]
Their lengths are at least \(2\eta\alpha_n/d_p\) and \(\xi\alpha_n/d_p\), respectively. Substituting \(\alpha_n=R_B^{-2n}/r_n\) and \(r_0\le r_n\le r_1\) proves the stated estimates with
\[
 c_0=\frac{\min\{P_-R_B^{-2},2\eta,\xi\}}{r_1},
 \qquad C_0=\frac{2T}{r_0}.
\]

For a sufficiently small \(\varepsilon>0\), set \(q_n=C_0R_B^{-2n}/d_p\) and choose \(n\) with \(q_n<\varepsilon\le q_{n-1}\). Both intervals lie in \((E_*,E_*+\varepsilon)\), and each has length at least \((c_0/C_0)q_n\ge \kappa\varepsilon\), where \(\kappa=c_0/(C_0R_B^2)\). The threshold for \(\varepsilon\) may depend on the background.
\end{proof}

\begin{proof}[\normalfont\bfseries Proof of Theorem~\ref{thm:main}]
For each fixed \(p\), Proposition~\ref{prop:routes} recovers the critical connection at energy zero for the active matrices and verifies the strict interval inclusions for every inactive matrix. By Lemma~\ref{pre:spectral}, the spectrum almost surely equals the same set for every \(0<\theta<1\). By Proposition~\ref{prop:transfer}, this spectrum contains zero and has infinitely many distinct bounded gaps accumulating at zero from the right.

For \(p\ge5\), \(L\) is the unique maximum in the period block. Any positive period \(d\) must preserve the residue class of this maximum under translation, so \(p\mid d\). At period four the same argument applies to the unique minimum \(\bar v_1-1\). At period three the entries are pairwise distinct. Finally, \eqref{eq:sparse} shows that each even period \(p\ge6\) has exactly six sites with nonzero background, and each odd period \(p\ge7\) has exactly seven. In both cases the background values belong to the set specified in the theorem.
\end{proof}

\subsection{Gap widths and locations at consecutive scales}

We apply Proposition~\ref{prop:transfer} at \(E_*=0\) to the backgrounds in Theorem~\ref{thm:main}, with \(R_B\) as in the proposition.

\begin{corollary}\label{cor:intro-quantitative}
Set \(s_{p,n}=R_B^{-2n}/p\). There are constants \(c_-,c_+,\kappa>0\), independent of \(p\), such that for each \(p\ge3\) there exist \(N(p)\in\mathbb N\), \(\varepsilon_0(p)>0\), and pairwise distinct bounded gaps \(G_{p,n}\subset(0,\infty)\) of \(\Sigma_p\) for all integers \(n\ge N(p)\). The following statements hold.
\begin{enumerate}[label=\textup{(\alph*)},leftmargin=*]
\item For \(n\ge N(p)\),
\[
\begin{aligned}
 c_-s_{p,n} &\le |G_{p,n}|\le c_+s_{p,n},\\
 c_-s_{p,n} &\le \inf G_{p,n}\le\sup G_{p,n}\le c_+s_{p,n}.
\end{aligned}
\]
\item For \(n>N(p)\), between \(G_{p,n}\) and \(G_{p,n-1}\) there is an interval in \(\Sigma_p\) of length at least \(c_-s_{p,n}\).
\item For \(0<\varepsilon<\varepsilon_0(p)\), the interval \((0,\varepsilon)\) contains an open interval of length at least \(\kappa\varepsilon\) disjoint from \(\Sigma_p\).
\end{enumerate}
\end{corollary}

\begin{proof}
By Proposition~\ref{prop:transfer} at \(E_*=0\), it remains to bound \(d_p/p\) above and below independently of \(p\).

For \(\ell=6,7\), take the solution \(\psi_\ell\) at the connection in \eqref{eq:transfer-wr}, where the periodic \(A\) tail on the right starts at site zero and \(\psi_\ell(0)=1\). In block \(i\), let \(j_i\) be the last interior site. Each new pair of sites with zero background is inserted immediately before \(j_i\). Set
\[
 c_\ell=\sum_{i\in\mathbb Z}
       \bigl(|\psi_\ell(j_i)|^2+|\psi_\ell(j_i-1)|^2\bigr).
\]
The solution decays exponentially on both periodic tails, so the series converges. Since two consecutive zero values would force the solution to vanish identically, \(c_\ell>0\). For the incoming state \((\psi_\ell(j_i),-\psi_\ell(j_i-1))^{\mathsf T}\), the transfer matrix at zero energy across an inserted pair with both noise variables zero is \(T_0^2=-I\). The amplitudes on the two inserted sites are \(\psi_\ell(j_i)\) and \(-\psi_\ell(j_i-1)\), up to a common sign. After \(k\) pairs have been inserted in each block, the absolute values of the amplitudes at the original sites are unchanged, with the same normalization at the start of the right tail. Equation~\eqref{eq:transfer-wr} now gives the exact identity
\begin{equation}\label{eq:width-derivative}
       d_{\ell+2k}=d_\ell+k c_\ell,\qquad \ell=6,7,\quad k\ge0 .
\end{equation}
Let \(D_-\) and \(D_+\) be the minimum and maximum of
\[
 \{d_j/j:3\le j\le7\}\cup\{c_6/2,c_7/2\}.
\]
These constants are positive. For \(p=\ell+2k\) with \(\ell=6,7\), formula~\eqref{eq:width-derivative} expresses \(d_p/p\) as the weighted average of \(d_\ell/\ell\) and \(c_\ell/2\), with weights \(\ell/p\) and \(2k/p\). Hence \(D_-p\le d_p\le D_+p\) for every \(p\ge3\).

Proposition~\ref{prop:transfer} gives the bounds for the gaps and the intervening spectral intervals with \(c_-=c_0/D_+\) and \(c_+=C_0/D_-\). Part (e) of the proposition gives the porosity bound at zero.
\end{proof}

\begin{remark}\label{rem:porosity-scale}
For fixed \(p\), the scales \(s_{p,n}\) of the selected gaps in \(\Sigma_p\) have common ratio \(R_B^{-2}\). This gives the uniform lower bound on right porosity. If porosity is measured using the radius of the omitted interval, the lower bound is \(\kappa/2\).
\end{remark}

In continuum periodization, Kirsch, Kotani and Simon obtain a gap-width lower bound of order \(1/L_j\) along a sequence of periods \(L_j\to\infty\), assuming nonzero reflection at a prescribed positive energy~\cite{KKS}. The geometric factor in Corollary~\ref{cor:intro-quantitative} describes the gaps within each fixed period.

\section{Four nonzero background sites per period via a resonant Dirichlet chain}\label{res:section}

\begin{proof}[\normalfont\bfseries Proof of Theorem~\ref{res:four}]
Use the parameters in~\eqref{eq:parameters} and let \(e\) be the connection energy of the model of period three. For \(p\ge5\), let \(m\) be the number of sites in the zero background chain and set
\[
 m=p-4,\qquad \alpha=\frac{\pi}{m+1},\qquad \gamma=-E_p=2\cos\alpha .
\]
Consider the block \(\widetilde V^{(p)}\) in~\eqref{eq:four-block}.

Since \(0\le\gamma<2<L\), the four entries outside the zero block have signs \(-,+,-,-\). The unique positive entry occurs once per block and forces minimal period \(p\).

For \(p=5\), we have \(E_5=0\) and \(\widetilde V^{(5)}\) is \(V^{(5)}\) in~\eqref{eq:sparse}, so Theorem~\ref{thm:main} gives the conclusion. Assume \(p\ge6\). At energy \(-\gamma\), the parameters at the last two sites agree with those of the model of period three at \(e\).

At energy \(-\gamma\), let \(x\) and \(y\) be the noiseless transfer parameters at the two ends of the prefix. Across these \(m+2\) sites the parameters are
\[
 x=e+L-\gamma,\quad \underbrace{-\gamma,\ldots,-\gamma}_{m},\quad y=-L .
\]
For \(k\ge-1\), set \(d_k=\sin((k+1)\alpha)/\sin\alpha\) for the recurrence on the free chain. The recurrence \(d_k=\gamma d_{k-1}-d_{k-2}\) gives
\[
 T_{-\gamma}^{m}
 =(-1)^m\begin{pmatrix}d_m&-d_{m-1}\\d_{m-1}&-d_{m-2}\end{pmatrix}
 =(-1)^m\begin{pmatrix}0&-1\\1&-\gamma\end{pmatrix}.
\]
If the noise vanishes on all \(m\) interior sites and the endpoint bits are \(i,j\in\{0,1\}\), then
\begin{equation}\label{res:identity}
 T_{y-\lambda j}T_{-\gamma}^{m}T_{x-\lambda i}
       =(-1)^m T_{e-\lambda(i+j)} .
\end{equation}
This follows by multiplication, since \(x+y+\gamma=e\). Thus for \(i=j=0\), the four choices on the last two sites recover the four active matrices of Lemma~\ref{lem:seed} up to the common sign.

Suppose now that at least one entry of the interior noise word \(\eta\in\{0,1\}^{m}\) is \(1\). Let \(\Delta_m\) denote the adjacency matrix of the path \(\{1,\ldots,m\}\), and set
\[
 K_\eta=\gamma I-\Delta_m+\lambda\operatorname{diag}(\eta).
\]
The unperturbed matrix \(\gamma I-\Delta_m\) is positive semidefinite with kernel spanned by \(s=(s_1,\ldots,s_m)^{\mathsf T}\), where \(s_k=\sin(k\alpha)\). Indeed, the vectors \((\sin(k\ell\alpha))_{k=1}^m\), \(1\le\ell\le m\), are eigenvectors of \(\Delta_m\) with distinct eigenvalues \(2\cos(\ell\alpha)\).

Choose \(j_0\) with \(\eta_{j_0}=1\). Let \(\delta_{j_0}\) be the \(j_0\)-th standard basis vector and set \(K^{(j_0)}=\gamma I-\Delta_m+\lambda\delta_{j_0}\delta_{j_0}^{\mathsf T}\). Since \(s_{j_0}>0\), the perturbation of rank one is strictly positive on the kernel. Hence \(K_\eta\ge K^{(j_0)}>0\).

The determinant of the unperturbed chain of length \(k\) is \(d_k=s_{k+1}/s_1\), with \(d_0=1\). A hat on a row or column index means that the row or column is deleted.

The determinant of the full chain is affine in its \(j_0\)-th diagonal entry. The coefficient of that entry is the product of the determinants of the two subchains on either side of \(j_0\). Using \(s_{m+1-j_0}=s_{j_0}\) and \(s_{m-j_0}=s_{j_0+1}\), we obtain
\begin{align*}
 \det K^{(j_0)}&=\lambda\,s_{j_0}^2/s_1^2,\\
 \det K^{(j_0)}_{\widehat m,\widehat m}
    &=1+\lambda\,s_{j_0}s_{j_0+1}/s_1^2 .
\end{align*}
The second determinant equals \(1\) when \(j_0=m\). Inversion reverses the order of positive definite matrices, so the cofactor formula yields
\begin{equation}\label{res:green}
 0<(K_\eta^{-1})_{mm}
 \le \frac{s_1^2}{\lambda s_{j_0}^2}
    +\frac{s_{j_0+1}}{s_{j_0}}
 \le \lambda^{-1}+\gamma .
\end{equation}
The last inequality uses \(s_{j_0}\ge s_1\) and \(s_{j_0+1}=\gamma s_{j_0}-s_{j_0-1}\), with \(s_0=0\).

For an incoming slope \(z\in[-3/2,104/1000]\), the slope after the left endpoint is \(-1/t\), where
\[
 t=x-\lambda i+z
 >\frac{568030}{10^6}+L-2-\lambda-\frac32>\frac25 .
\]
The slope after the \(k\)-th interior site satisfies \(u_k=(\gamma+\lambda\eta_k-u_{k-1})^{-1}\), with \(u_0=-1/t\). Taking Schur complements in the tridiagonal matrix gives
\[
 u_m=\big[(K_\eta+t^{-1}\delta_1\delta_1^{\mathsf T})^{-1}\big]_{mm}.
\]
Indeed, the lower right entries of the inverses of the successive leading principal submatrices obey the same recurrence. Each submatrix is positive definite, so every denominator is positive. Applying the order reversal under inversion once more and using~\eqref{res:green}, we obtain \(0<u_m\le\gamma+\lambda^{-1}\). The outgoing slope \(s_{\rm out}\) at the right endpoint satisfies
\[
 \frac{36}{1000}
 <\frac1{L+\lambda}\le
 s_{\rm out}=\frac1{L+\lambda j-u_m}
 <\frac1{L-2-\lambda^{-1}}<\frac{75}{1000}.
\]
If all interior bits vanish and \(i+j\ge1\), then~\eqref{res:identity} and the slope estimate in the proof of Lemma~\ref{lem:five} give \(s_{\rm out}\in(0.039,0.088)\).
Thus every prefix with a nonzero noise bit sends \(\mathcal I\) into \((0.036,0.09)\). At \(E=-\gamma\), \([w,H]\subset\mathcal I\), and Lemma~\ref{lem:tail} puts the full image compactly into one of the three target intervals.

The four words whose prefixes contain only zeros share the same prefix matrix at all energies. The factorization in Proposition~\ref{prop:routes} then gives \(RA^{-1}C=B\). Equation~\eqref{res:identity} identifies their matrices at \(-\gamma\) with the seed up to a common sign. Proposition~\ref{prop:transfer} gives \(-\gamma\) as a spectral point with distinct bounded gaps approaching it from the right.
\end{proof}

\begin{corollary}\label{cor:resonant-gaps}
Set \(\widetilde s_{p,n}=R_B^{-2n}/p^3\). For the background in Theorem~\ref{res:four}, use the same \(R_B\) as in Corollary~\ref{cor:intro-quantitative}. At each index \(n\), use the same block order \(B^nC^2A^{m_n}\) and the same counts \((n,2,m_n)\), with the blocks defined by the background of Theorem~\ref{res:four}. There are constants \(c_-,c_+>0\), independent of \(p\), such that for each \(p\ge5\) there exist an integer \(\widetilde N(p)\) and pairwise distinct bounded gaps \(\widetilde G_{p,n}\subset(E_p,\infty)\) of \(\widetilde\Sigma_p\) for all integers \(n\ge\widetilde N(p)\). The following statements hold.
\begin{enumerate}[label=\textup{(\alph*)},leftmargin=*]
\item For \(n\ge\widetilde N(p)\),
\[
\begin{aligned}
 c_-\widetilde s_{p,n} &\le |\widetilde G_{p,n}|\le c_+\widetilde s_{p,n},\\
 c_-\widetilde s_{p,n} &\le \inf\widetilde G_{p,n}-E_p
 \le \sup\widetilde G_{p,n}-E_p\le c_+\widetilde s_{p,n}.
\end{aligned}
\]
\item For \(n>\widetilde N(p)\), between \(\widetilde G_{p,n}\) and \(\widetilde G_{p,n-1}\) there is an interval in \(\widetilde\Sigma_p\) of length at least \(c_-\widetilde s_{p,n}\).
\end{enumerate}
\end{corollary}

\begin{proof}
For \(p=5\), Corollary~\ref{cor:intro-quantitative} gives the stated bounds after rescaling the constants by \(25\), since \(5^{-1}=25\cdot5^{-3}\). For \(p\ge6\), use Proposition~\ref{prop:transfer} at \(E_*=E_p\) as in the preceding proof. Retain \(m,\alpha,\gamma,x\) from that proof and put \(b=v_1\). Here \(d_p\) again denotes the energy derivative in~\eqref{eq:transfer-wr}.

Take the connection solution \(\psi\) of~\eqref{eq:transfer-wr} at \(E_p\), with \(\psi(0)=1\) at the start of the first \(A\) block. Set \(X_i=(a_i,b_i)^{\mathsf T}=(\psi(ip),-\psi(ip-1))^{\mathsf T}\). At each block entrance the state agrees up to sign with that of the normalized connection for period three, by~\eqref{res:identity}. The two periodic tails decay exponentially in the block index, so \(\sum_i\|X_i\|^2\) is finite and independent of \(p\).

In block \(i\), set \(w_j=(-1)^j\psi(ip+j)\) and \(q_i=(e+L-\cos\alpha)a_i+b_i\). The first step gives \(w_0=a_i\), \(w_1=-xa_i-b_i\). On the chain, \(w_{j+1}=\gamma w_j-w_{j-1}\), hence
\[
 w_j=a_i\cos(j\alpha)-q_i\,\frac{\sin(j\alpha)}{\sin\alpha}
 \qquad (0\le j\le m+1).
\]
The mixed term vanishes upon pairing \(j\) with \(m+1-j\). Since \(\sum_{j=1}^m\sin^2(j\alpha)=(m+1)/2\), the mass on the chain is
\begin{equation}\label{eq:resonant-mass}
 \sum_{j=1}^m|\psi(ip+j)|^2
 =\frac{m-1}{2}a_i^2+\frac{m+1}{2\sin^2\alpha}q_i^2 .
\end{equation}
At the four other sites of the block the amplitudes are bounded by a constant times \(\|X_i\|\), uniformly in \(p\). Indeed, \(w_{m+1}=-a_i\), the state at site \(ip+m+2\) is \((-1)^mT_eX_i\), and the next step uses one of the two fixed parameters \(e-b\) and \(e-b-\lambda\).

Since \(m+1=p-3\) and \(\sin\alpha\asymp p^{-1}\), summing over all blocks gives \(\|\psi\|^2\le C p^3\). For the reverse bound it suffices to take block zero. There \((a_0,b_0)=(1,s_A)\) and \(q_0=e+L-\cos\alpha+s_A>14\) by~\eqref{eq:seed-radii}. Thus \(d_p=\|\psi\|^2\asymp p^3\), with constants independent of \(p\). Proposition~\ref{prop:transfer} now gives the stated estimates for the gaps and the intervening spectral intervals.
\end{proof}

\begin{remark}\label{rem:dirichlet-resonance}
Both families use the geometric scale \(R_B^{-2n}\). By Proposition~\ref{prop:transfer}, their energy scales are comparable to \(R_B^{-2n}/d_p\), where \(d_p\) is the energy derivative of the difference of the connecting slopes at the accumulation energy. Equation~\eqref{eq:transfer-wr} identifies this derivative with the squared norm of the normalized connection solution. At zero, each inserted pair adds a fixed positive amount to that squared norm within each parity family, giving \(d_p\asymp p\) by~\eqref{eq:width-derivative}. At \(E_p\), resonance gives a sine term of amplitude of order \(p\) on the chain in the first block of the right periodic tail. The sum of squares along the chain gives \(d_p\asymp p^3\) by~\eqref{eq:resonant-mass}.

For fixed \(p\ge5\) and \(n\ge\max\{N(p),\widetilde N(p)\}\), the two constructions use the same numbers of blocks in the periodic word \(B^nC^2A^{m_n}\). Corollaries~\ref{cor:intro-quantitative} and~\ref{cor:resonant-gaps} therefore give
\[
 \frac{|\widetilde G_{p,n}|}{|G_{p,n}|}\asymp p^{-2},\qquad
 \frac{\inf\widetilde G_{p,n}-E_p}{\inf G_{p,n}}\asymp p^{-2},
\]
with comparison constants independent of \(p\) and \(n\).
\end{remark}

\appendix

\section{Estimates for the integer scales}\label{app:width}

In this appendix, \(R_A,R_B\) are the expanding eigenvalue moduli at the connection energy for the model of period three. The constants \(P_\pm,Q_\pm\) are the coefficients in the projective expansions~\eqref{eq:transfer-asym}, and \(a^*,b^*\) are the corresponding trace coefficients in~\eqref{eq:transfer-trace}.

\begin{lemma}\label{lem:width-parameters}
Take
\[
 r_0=0.052,\qquad r_1=0.111,\qquad \eta=0.000004,\qquad T=0.2 .
\]
Then
\[
 \begin{gathered}
  R_A^2r_0<r_1,\qquad R_B>10.4336,\qquad P_->0.018,\qquad a^*+r_1 b^*<T,\\
  \min\{(Q_+-Q_-)r_0-P_-,\ P_+-P_--Q_-r_1\}>3\eta .
 \end{gathered}
\]
\end{lemma}

\begin{proof}
Set \(a=u_A^*,\,s=s_A^*,\,u=u_B^*,\,S=-s_B^*,\,\Delta=s-a\), and let \(z_-=C^*(-1/4)\) and \(z_+=C^*0\). Table~\ref{tab:seed-bounds} gives the rational bounds at the connection energy. In~\cite{CL}, \(t\) denotes our \(z_+\). For \(t=-307/200\) used here, Lemma~\ref{lem:seed} gives \(0<B^*t<0.08734\). Here \(h_M(v)=M_{11}+M_{12}v\), and we omit the stars on matrix entries.

\begin{table}[htbp]
\centering
\small
\renewcommand{\arraystretch}{1.14}
\begin{tabularx}{\linewidth}{@{}l>{\raggedright\arraybackslash}X@{\hspace{1.8em}}l>{\raggedright\arraybackslash}X@{}}
\toprule
Quantity & Rational interval & Quantity & Rational interval \\
\midrule
\(a\) & \((-0.55183,-0.55181)\) & \(s\) & \((-0.51438,-0.51436)\) \\
\(u\) & \((0.10258,0.10260)\) & \(R_A^2\) & \((1,2.13)\) \\
\(z_-\) & \((-0.47625,-0.47623)\) & \(z_+\) & \((-0.42954,-0.42953)\) \\
\(B_{11}\) & \((10.42449,10.42450)\) & \(B_{12}\) & \((0.08894,0.08896)\) \\
\(-C_{11}\) & \((2.466911,2.466914)\) & \(-C_{12}\) & \((1.188107,1.188110)\) \\
\(B_{22}=C_{22}\) & \([0.104969,0.104970]\) &  &  \\
\bottomrule
\end{tabularx}
\caption{Rational bounds at \(E=e\) from~\cite{CL}}\label{tab:seed-bounds}
\end{table}

The sum of the roots of the equation for the fixed points of \(B\) gives \(S=u+(B_{11}-B_{22})/B_{12}\). Hence \(116<S<117\), and \(R_B=h_B(u)>10.4336\). Since \(D=C^2=(\tr C)C-I\),
\[
 h_D(u)=(-C_{11}-C_{22})(-C_{11}-C_{12}u)-1 .
\]
Both factors are positive, so the entry bounds give
\[
\begin{aligned}
 5.1145
 &<(2.466911-0.104970)(2.466911+1.188107\cdot 0.10258)-1\\
 &<h_D(u)<(2.466914-0.104969)(2.466914+1.188110\cdot 0.10260)-1
 <5.1147 .
\end{aligned}
\]
It follows that \(d_0=(D^*)'(u)=h_D(u)^{-2}\) satisfies \(0.038225<d_0<0.03823\).

The formula \(P_\pm=\Delta(z_\pm-s)/(z_\pm-a)\) yields
\[
\begin{aligned}
 0.018 &<\frac{0.03743\cdot 0.03811}{0.07560}
       <P_-\le \frac{z_--a}{4}<0.01890,\\
 P_+-P_-&=\frac{\Delta^2(z_+-z_-)}{(z_+-a)(z_--a)}
     >\frac{0.03743^2\cdot 0.04669}{0.12230\cdot 0.07560}>0.00707 .
\end{aligned}
\]
The upper bound in the first line uses \(4(s-a)(z_--s)\le (z_--a)^2\). For the other two coefficients we have
\[
 \frac{Q_-}{d_0}=(u-t)\left(1+\frac{u-t}{S+t}\right),\qquad
 \frac{Q_+}{d_0}=R_B^2u\left(1+\frac{u}{S}\right).
\]
Substituting the bounds on \(S\) and \(u\) gives
\[
\begin{aligned}
 Q_- &<0.03823\cdot1.63760\left(1+\frac{1.63760}{114.465}\right)<0.06352,\\
 \frac{Q_+-Q_-}{d_0}
 &>10.4336^2\cdot 0.10258\left(1+\frac{0.10258}{117}\right)
      -1.63760\left(1+\frac{1.63760}{114.465}\right)>9.515 .
\end{aligned}
\]
The two margins in the scale window satisfy
\[
\begin{aligned}
 (Q_+-Q_-)r_0-P_-&>0.038225\cdot9.515\cdot0.052-0.01890=0.0000129655,\\
 P_+-P_--Q_-r_1&>0.00707-0.06352\cdot0.111=0.00001928 .
\end{aligned}
\]
Both exceed \(3\eta=0.000012\). Also, \(R_A^2r_0<2.13\cdot0.052=0.11076<r_1\).

From \(0.03743<\Delta<0.038\), \(S+a>115\), and the bound on \(d_0\), the trace coefficients satisfy
\[
\begin{gathered}
 P_+<a^*=\Delta+\frac{\Delta^2}{S+a}<0.039,\\
 0<b^*=d_0(u-a)\frac{S+u}{S+a}
      <0.03823\cdot 0.655\left(1+\frac{0.655}{115}\right)<0.026 .
\end{gathered}
\]
Hence \(a^*+r_1 b^*<0.041886<T\). For the adjacent scales in~\eqref{eq:transfer-integers}, \(r_1/(r_0 R_B^2)<0.01961<0.02\) and \(2T\cdot0.02=0.008<P_-\).
\end{proof}

\section{Algebraicity of the connection energy}\label{app:algebraicity}

The connection energy \(e\) in Lemma~\ref{lem:seed} is a real algebraic number.

\begin{proof}
With the matrices of Section~\ref{sec:seed}, set \(U=DBD^{-1}\) and \(P(E)=\operatorname{tr}(AUA^{-1}U^{-1})-2\). All entries in the commutator lie in \(\mathbb Z[1/10][E]\), since the entries of \(A,B,C\) lie there and the matrices have determinant one. Here \(\mathbb Z[1/10]=\{a/10^k:a\in\mathbb Z,\ k\in\mathbb Z_{\ge0}\}\).

At \(E=e\), the common eigenline of \(A\) and \(U\) puts the commutator in upper triangular form with both diagonal entries equal to one. Thus \(P(e)=0\). To check that \(P\) is nonzero, evaluate at \(E=0\) and reduce modulo \(3\), using \(10\equiv1\pmod3\). Formula~\eqref{eq:seed-product} gives
\[
 A=\begin{pmatrix}0&-1\\1&0\end{pmatrix},\qquad
 U=\begin{pmatrix}1&-1\\0&1\end{pmatrix},\qquad
 AUA^{-1}U^{-1}=\begin{pmatrix}1&1\\1&2\end{pmatrix}
 \quad\text{in }\mathbb Z_3=\mathbb Z/3\mathbb Z .
\]
Thus \(P(0)\equiv1\pmod3\), so \(P\) is nonzero and \(e\) is algebraic.
\end{proof}

\end{document}